\documentclass[11pt, reqno]{amsart}

\usepackage{amssymb,geometry,stmaryrd,color}
\usepackage[T1]{fontenc}

\makeatletter
\@namedef{subjclassname@2020}{\textup{2020} Mathematics Subject Classification}
\makeatother

\evensidemargin\oddsidemargin

\begin{document}

\baselineskip=18pt
\setcounter{page}{1}
    
\newtheorem{Conj}{Conjecture}
\newtheorem{TheoA}{Theorem A\!\!}
\newtheorem{TheoB}{Theorem B\!\!}
\newtheorem{TheoC}{Theorem C\!\!}
\newtheorem{TheoD}{Theorem D\!\!}
\newtheorem{Lemm}{Lemma}
\newtheorem{Rem}{Remark}
\newtheorem{Def}{Definition}
\newtheorem{Coro}{Corollary}
\newtheorem{Propo}{Proposition}

\renewcommand{\theConj}{}
\renewcommand{\theTheoA}{}
\renewcommand{\theTheoB}{}
\renewcommand{\theTheoC}{}
\renewcommand{\theTheoD}{}

\def\a{\alpha}
\def\b{\beta}
\def\g{\gamma}
\def\B{{\bf B}} 
\def\CC{{\mathbb{C}}} 
\def\cG{{\mathcal{G}}} 
\def\cB{{\mathcal{B}}} 
\def\cI{{\mathcal{I}}} 
\def\cS{{\mathcal{S}}}
\def\UU{{\mathcal{U}}}
\def\ca{c_{\a}}
\def\ka{\kappa_{\a}}
\def\coa{c_{\a, 0}}
\def\cua{c_{\a, u}}
\def\cL{{\mathcal{L}}} 
\def\cW{{\mathcal{W}}} 
\def\Ea{E_\a}
\def\Eab{E_{\a,\b}}
\def\eps{{\varepsilon}} 
\def\esp{{\mathbb{E}}} 
\def\Ga{{\Gamma}} 
\def\G{{\bf G}} 
\def\K{{\bf K}}
\def\HH{{\bf H}}
\def\ii{{\rm i}}
\def\e{{\rm e}}
\def\L{{\bf L}}
\def\lbd{\lambda}
\def\lacc{\left\{}
\def\lcr{\left[}
\def\lpa{\left(}
\def\lva{\left|}
\def\M{{\bf M}}
\def\T{{\bf T}}
\def\Ma{\M_\a}
\def\Mab{\M_{\a,\b}}
\def\NN{{\mathbb{N}}} 
\def\pb{{\mathbb{P}}}
\def\pa{{\varphi_\a}}
\def\paa{{\varphi_{\a,1-\a}}} 
\def\pab{{\varphi_{\a,\b}}} 
\def\tpab{\hat{\varphi}_{a,b}} 
\def\tpa{\hat{\psi}_{\a}}
\def\tppa{\tilde{\psi}_{\a}} 
\def\tva{\hat{\varphi}_{\a}} 
\def\rl{{\mathbb{R}}}
\def\racc{\right\}}
\def\rpa{\right)}
\def\rcr{\right]}
\def\rva{\right|}
\def\prost{{\succ_{\! st}}}
\def\W{{\bf W}}
\def\X{{\bf X}}
\def\Z{{\bf Z}}
\def\Xab{\X_{\a,\b}}
\def\XX{{\mathcal X}}
\def\Y{{\bf Y}}
\def\U{{\bf U}}
\def\V{{\bf V}}
\def\Un{{\bf 1}}
\def\ZZ{{\mathbb{Z}}}
\def\A{{\bf A}}
\def\AA{{\mathcal A}}
\def\hAA{{\hat \AA}}
\def\hL{{\hat L}}
\def\hT{{\hat T}}

\def\claw{\stackrel{d}{\longrightarrow}}
\def\elaw{\stackrel{d}{=}}
\def\qed{\hfill$\square$}

\newcommand*\pFqskip{8mu}
\catcode`,\active
\newcommand*\pFq{\begingroup
        \catcode`\,\active
        \def ,{\mskip\pFqskip\relax}%
        \dopFq
}
\catcode`\,12
\def\dopFq#1#2#3#4#5{%
        {}_{#1}F_{#2}\biggl[\genfrac..{0pt}{}{#3}{#4};#5\biggr]%
        \endgroup
}

\title{Mittag-Leffler functions and convex ordering}

\author[Rui. A. C. Ferreira]{Rui A. C. Ferreira}

\address{CMUP, Faculdade de Ciências, Universidade do Porto,
Rua do Campo Alegre 687, 4169--007 Porto, Portugal.  {\em Email}: {\tt rui.ac.ferreira@fc.up.pt}.}

\author[Thomas Simon]{Thomas Simon}

\address{Laboratoire Paul Painlev\'e, UMR 8524, Universit\'e de Lille, 42 rue Paul Duez, 59000 Lille, France. {\em Email}: {\tt thomas.simon@univ-lille.fr}}



\begin{abstract} 
The monotonicity of the Mittag-Leffler function $E_{\alpha}$ with respect to the parameter $\alpha$ is investigated, via some convex ordering properties for related random variables. In particular, it is shown that the mapping $\alpha\mapsto E_\alpha(x^\alpha)$ decreases on $(0,2)$ for all $x> 0$, that the mapping $\alpha\mapsto E_\alpha(-x^\alpha)$ decreases on $(0,1)$ for all $x\ge 1$ and that the mapping $\alpha\mapsto E_\alpha(\Gamma(1+\alpha)x)$ decreases on $(0,1)$ for all $x\in{\mathbb R}^\ast.$ Analogous results are presented for the two parameter Mittag-Leffler functions $E_{\alpha, \beta}$ with $\beta\ge \alpha,$ with an emphasis on the extremal case $\beta =\alpha.$ Several applications of these results are discussed for Abelian integral equations and subdiffusions.
\end{abstract}

\maketitle

\section{Introduction and presentation of the results}

Introduced in \cite{ML}, the classical Mittag-Leffler function
$$\Ea(x)\, =\, \sum_{n\ge 0} \frac{x^n}{\Ga(1+ n\a)}$$
with parameter $\a > 0,$ is one of the simplest special functions. In the case when $\a$ is rational, this function has a hypergeometric character and solves a certain linear ordinary differential equation. In general, the Mittag-Leffler function solves an Abelian integral equation with a fractional kernel and for this reason, it has been the object of many studies in the literature on fractional calculus. See the recent monographs \cite{AMP, GKMR} for comprehensive accounts. Meanwhile, some basic properties of these functions are still unsettled. 

In this paper, we wish to study some monotonicity properties of the Mittag-Leffler function with respect to the fractional parameter $\a.$ This problem is classical for the usual special functions, in relationship with various comparison inequalities for the involved differential equations. As an instance, let us refer to our previous paper \cite{FS1} for such comparison inequalities in the framework of confluent hypergeometric functions. For the Mittag-Leffler functions, the second author had obtained in \cite{Frechet} various stochastic ordering results for related random variables, leading to the hyperbolic bounds 
\begin{equation}
\label{Unif}
\frac{1}{1+ \Ga(1-\a) x}\, \le \, \Ea(-x)\, \le\, \frac{1}{1 + \frac{x}{\Ga(1+\a)}}
\end{equation}
for all $\a\in (0,1)$ and $x\ge 0.$ These bounds, which were later on generalized in \cite{BS} for the two parameter Mittag-Leffler functions and the more general Kilbas-Saigo functions - can then be used for several problems such as Van der Corput estimates \cite{RT}, inverse problems for subdiffusions \cite{JZ} or boundedness properties for fractional differential equations \cite{DNNT}.

The purpose of this paper is to apply the notion of convex ordering for studying the monotonicity of the mapping $\a\mapsto\Ea,$ in various contexts. Let us first recall that for two real random variables $X$ and $Y,$ the convex ordering relationship
$$X\,\prec_{cx}\, Y$$
means that $\esp[\varphi(X)]\le\esp[\varphi(Y)]$ for all convex functions such that the two expectations exist. By the convexity of the exponential function, this property has an immediate application for comparing Laplace transforms or moment generating functions, which is a known feature of several Mittag-Leffler functions or transformations thereof - see \cite{MLCM} for an account. The methodology we use to obtain convex orderings in this paper is classical, and makes an extensive use of single and double intersection properties, as displayed in the standard reference \cite{Shaked} on the topic. Considering the sometimes technical arguments and the specificity of some statements in \cite{Frechet}, our proofs will be more transparent and our findings will also have a simpler character. Our first result is the following. 

\begin{TheoA}
For every $x > 0,$ the mapping $\alpha\mapsto E_\alpha(x^\alpha)$ decreases on $(0,2)$.
\end{TheoA}

The proof of the above result is different according as $\a\in (0,1)$ or $\a\in (1,2).$ In case $\a\in (0,1),$ the function $\Ea(x^\a)$ is the difference of two Laplace transforms and we can apply a complete monotonicity result of \cite{Alzer}, which basically amounts to a non-intersection property for the kernels in one Laplace transform. In case $\a\in(1,2),$ the function $\Ea(x^\a)$ is the moment generating function of a centered random variable with an atom, and we show a double intersection property with the help of the generalized Descartes' rule of signs. If we consider the two-parameter Mittag-Leffler function 
$$\Eab (x)\, =\, \sum_{n\ge 0} \frac{x^n}{\Ga(\b +\a n)}$$
with $\a,\b > 0,$ a direct consequence of Theorem A and of the interpolation formula
$$\Eab(x^\a)\, =\, \frac{1}{\Ga(\b -1)} \int_0^1 (1-t)^{\b -2} \Ea((tx)^\a)\, dt$$
which is valid for all $\b > 1$ and $x \ge 0,$ is that the mapping $\a\mapsto\Eab(x^\a)$ also decreases on $(0,2)$ for all $\b > 1$ and $x >0.$ However, we shall observe in Remark \ref{Rem0} below that the mapping $\a\mapsto \Eab(x^\a)$ may not be monotonic for $\b < 1$, and also that the mapping $\a\mapsto\Ea(x^\a)$ has increase points for large $\a$ and some specific values of $x.$ In case $\a\in (0,1),$ the derivative of the function $x\mapsto\Ea(\lbd x^\a)$ appears as the fundamental solution to Abelian integral equations of the second kind, and we show in Corollaries \ref{FDE1} and \ref{FDE2} that the proof of Theorem A leads to interesting comparison theorems for such equations, in the expanding case $\lbd > 0.$ 

\medskip

In the literature on fractional calculus, the completely monotonic function $x\mapsto \Ea(-x^\a)$ with $\a\in (0,1)$ plays a central role and we refer again to \cite{GKMR} for an account - see in particular Chapter 8.1 therein for its connections with fractional relaxation. Our next main result shows that monotonicity with respect to the fractional parameter has a very different character for such functions.  

\begin{TheoB}
For every $0 <\a <\b <1,$ there exists $x_{\a,\b} \in (0,1)$ such that the function 
$$x\,\mapsto\,\Ea(-x^\a) - E_\b(-x^\b)$$ 
is negative for $x\in (0,x_{\a,\b})$ and positive for $x\in (x_{\a,\b},\infty).$ Moreover, one has the bounds
$$ \lpa\frac{1}{\Ga(1-\b)\Ga(1+\a)}\rpa^{\frac{1}{\b-\a}} <\; x_{\a,\b} \; <\,  \lpa\frac{\Ga(1+\b)}{\Ga(1+\a)}\rpa^{\frac{1}{\b-\a}}\!.$$
\end{TheoB}

This result implies that for every $x\ge 1,$ the mapping $\alpha\mapsto E_\alpha(-x^\alpha)$ decreases on $(0,1)$. This completes the main result of the recent paper \cite{AS} - see Theorem 3.5 therein, which shows that for every $\a_\ast\in (0,1)$ there exists $c_\ast \in (0,1)$ such that the mapping $\alpha\mapsto E_\alpha(-x^\alpha)$ increases on $(\a_\ast, 1)$ for every $x\in (0, c_\ast)$. In \cite{AS}, this property is applied to the inverse problem of recovering the fractional exponent of subdiffusions from observed data - see Theorem 5.4 therein. We believe that our uniform monotonicity result could be useful for such questions, which we have not however investigated as yet. 

Another consequence of Theorem B is an alternative proof of Theorem A in the case $\a\in (1,2)$ - see Remark \ref{Rem1}, and a precise reason why the main result of \cite{Alzer} must fail in general - see Proposition \ref{SC2}. We also show in Proposition \ref{Eaa1} that the single intersection property of Theorem B can be obtained from the double intersection property for the derivative $x^{\a-1}E_{\a,\a} (-x^\a),$ which however does not imply the important bound $x_{\a,\b} < 1.$ Finally, we observe in Proposition \ref{Bound} that the material used for Theorem B also leads for $\a\in (0,4]$ to optimal bounds for the difference $\a\Ea(x) -e^{x^{1/\a}}$ between the Mittag-Leffler function and its leading term in the expansion at $\infty.$  

\medskip

Our third main result is the following improvement on the bounds given in \eqref{Unif}. 

\begin{TheoC}
For every $x\in\rl^\ast,$ the function $\alpha\mapsto E_\alpha\lpa \Ga(1+\a) x\rpa $ decreases on $(0,1)$.
\end{TheoC}

Applied on the negative half-line, this result yields the uniform bounds
\begin{equation}
\label{Unif1}
e^{-x}\, \le\, E_\b\lpa - \Ga(1+\b) x\rpa\, \le\, E_\a\lpa-\Ga(1+\a) x\rpa\, \le\, \frac{1}{1+x}
\end{equation}
for all $0 < \a < \b < 1$ and $x\ge 0,$ which shows that the optimal upper bound in \eqref{Unif} can be obtained as an increasing limit when $\a\downarrow 0,$ and also that the lower bound in \eqref{Unif}, which is optimal when $x\to \infty,$ can be improved when $x\to 0$. These bounds also specify how the function $\Ea\lpa -\Ga(1+\a) x^2\rpa$ interpolates between the Lorentzian curve $1/(1+x^2)$ at the limit $\a \to 0$ and the Gaussian curve $e^{-x^2}$ at the limit $\a\to 1.$ This interpolation turns out to be monotonic after a rescaling that equalizes all second derivatives at zero. Observe that the animation on this interpolation for the function $\Ea(-x^2)$ which is given at the beginning of \cite{Wiki} is not, despite appearances, monotonic.  

Theorem C is obtained as a consequence of the convex ordering for the rescaled Mittag-Leffler random variables $\Ga(1+\a)\M_\a,$ whose moment generating function is $E_\alpha\lpa \Ga(1+\a) x\rpa$ - see the introduction of our previous paper \cite{FS2} and the references therein for more details on these random variables. This convex ordering had already been derived on $(1/2,1)$ in \cite{Frechet} via some discrete factorization with Beta and Gamma random variables. In the present paper, we use a factorization with Beta random variables only and an exact power transform which were already applied in \cite{Sun} for other purposes, and turn out to be particularly adapted to our problem.   

The latter power beta factorization extends for all $\a\in (0,1)$ and $\b\ge\a$ to the generalized Mittag-Leffler random variables $\M_{\a,\b}$ having moment generating function $\Eab(x)$  - see again the introduction of \cite{FS2} for more details on such random variables. And it is equally useful for proving our last main result dealing with the monotonicity of $\Eab$ with respect to the second parameter $\b$. 

\begin{TheoD}
\label{Eab}
For every $\a\in (0,1]$ and $x\in\rl^\ast,$ the mapping 
$$\b\,\mapsto\, \Ga(\b)\, E_{\a,\b} \lpa \frac{\Ga(\b+\a)}{\Ga(\b)}\, x\rpa$$
increases on $(\a,\infty).$ Moreover, it converges to $\infty$ if $x\ge 1$ and to $1/(1-x)$ if $x < 1.$ 
\end{TheoD}

This result is also obtained via convex ordering, in a different manner than Theorem C: the argument involves here the factorization of $\M_{\a,\b}$ with the "extremal" distribution $\M_{\a,\a}$ and the power transformation of a single beta distribution. Similarly as above, Theorem D improves on the hyperbolic bounds obtained in \cite{BS} for the two parameter Mittag-Leffler functions - see \eqref{Bd1} and Remark \ref{Rem4} below. In Proposition \ref{Eaa}, we also study the monotonicity with respect to the parameter $\a$ of the suitably rescaled extremal Mittag-Leffler function $E_{\a,\a}$ on the positive half-line, which leads to an improvement on the Cauchy bounds for the latter function - see again Remark \ref{Rem4}. Finally, in Proposition \ref{Inv} we show how an argument during the proof of Theorem B combined with these optimal Cauchy bounds can be applied for a precise variation analysis of the function $x\mapsto xE_{\a,\a} (-x),$ which is crucial in the sensitivity analysis for subdiffusions in the recent paper \cite{JZ}.

\section{Proof of Theorem A}
\label{A}

We start with the case $\alpha\in(0,1),$ which is a consequence of the following integral representation:
\begin{equation}\label{M-L Repr}
    E_\alpha(x^\alpha)\,=\, \frac{e^x}{\alpha}\,-\,\int_0^\infty e^{-xt}f_\a(t)\,dt, \qquad x\ge 0,
    \end{equation}
with
$$f_\a(t)\,=\, \frac{\sin(\pi\alpha)\, t^{\alpha-1}}{\pi(t^{2\alpha}-2\cos(\pi\alpha)t^\alpha+1)}$$ 
for all $t > 0.$ This representation is given in Formula (3.4) of \cite{MLCM} - see also Formulas (1.1.12) and (1.1.16) with $\rho =1/\alpha$ and $\mu = 1$ in \cite{Popov}. If $0<\alpha<\beta < 1,$ it follows from Lemma 2.2 and Formula (3.3) in \cite{Alzer} that 
$$f_{\alpha}(t)-f_{\beta}(t)\, > \, \left(\frac{\beta}{\alpha}-1\right)f_{\beta}(t)\,> \,0, \qquad t > 0.$$
Therefore, we obtain
\begin{eqnarray*}
        E_{\alpha}(x^{\alpha})\,-\,E_{\beta}(x^{\beta})& = & \frac{e^x}{\alpha}-\frac{e^x}{\beta}\, +\, \int_0^\infty e^{-xt}\lpa f_{\beta}(t)-f_{\alpha}(t)\rpa \,dt\\
        & > &  \frac{e^x}{\alpha}-\frac{e^x}{\beta}\, +\, \int_0^\infty \lpa f_{\beta}(t)-f_{\alpha} (t)\rpa\, dt\; = \; \frac{(e^x-1)(\beta-\alpha)}{\alpha\beta}\; > \; 0
\end{eqnarray*}
as required, where for the second equality we have used \eqref{M-L Repr} with $x = 0$. 

We next proceed to the case $\alpha\in (1,2),$ which relies on convex ordering. Formula (3.6) in \cite{MLCM} shows that \eqref{M-L Repr} remains true in this context, but the function $f_\a$ becomes here everywhere negative so that the formula has a different meaning. Specifically, if we introduce the real random variable $\X_\a$ whose law is given by
$$\frac{1}{\alpha}\, \delta_{-1} (dt)\, +\, \lpa\frac{(-\sin(\pi\a)) t^{\a -1}}{\pi (t^{2\a} - 2\cos(\pi\a) t^\a +1)}\, \Un_{\{t > 0\}}\rpa dt,$$
then \eqref{M-L Repr} reads as
\begin{equation}\label{M-L Rep}
E_\alpha(x^\alpha)\, =\, \esp[e^{-x\X_\a}], \qquad x\ge 0.
\end{equation}
Observe that $\X_\a$ is integrable with $\esp[\X_\a] = 0$ for all $\a\in (1,2)$ because
$$\frac{d}{dx} \lpa E_\alpha(x^\alpha)\rpa_{x=0} \, =\, 0.$$
Since $y\mapsto e^{-xy}$ is strictly convex for every $x > 0,$ it suffices to show that $\X_{\beta}\prec_{cx} \X_{\alpha}$
if $1 < \alpha <\beta <2.$ Setting $F_\a(t) = \pb[\X_\a \le t]$ and $F_\beta(t) = \pb[\X_\b \le t]$ for all $t\in\rl,$ it is enough to prove that the function
$$t\,\mapsto\, F_\a(t)\, -\, F_\b(t)$$
vanishes once on the common support $[-1,\infty),$ starting positive. To see this, apply Theorem 3.A.44 in \cite{Shaked} or use alternatively the formula
$$\frac{1}{x}\lpa \Ea(x^\a)\, -\, E_\b(x^\b)\rpa\, =\, \int_{-1}^\infty e^{-xt} \lpa F_\a(t)\, -\, F_\b(t)\rpa\, dt$$ combined with Theorem 3.1.(b) p.21 in \cite{Karlin} and the fact that the function on the left-hand side vanishes at zero and $\infty.$ The function $F_\a(t)-F_\b(t)$ is clearly positive on $[-1,0],$ smooth on $(0, \infty)$ and converges to zero as $x\to\infty.$ Hence, by an elementary analysis - see also the proof of Theorem 3.A.44 in \cite{Shaked}, we are reduced to show that the function
$$F_{\a,\b}(t) \, =\, F_\a'(t)\, -\, F_\beta'(t)$$
vanishes twice on $(0,\infty),$ starting positive. Computing
$$F_{\a,\b}(t)\, =\, \frac{(-\sin(\pi\alpha))\, t^{\alpha}(1+ t^{2\beta}) \, + \, \sin(\pi\beta)\, t^{\beta}(1+ t^{2\a}) \, +\, 2\sin(\pi(\a-\beta))\, t^{\a+\beta}}{\pi \,t(t^{2\a} - 2\cos(\pi\a) t^\a +1)(t^{2\beta} - 2\cos(\pi\beta) t^\beta +1)}$$
shows that the number of zeroes of $F_{\a,\b}$ on $(0,\infty)$ is that of the function
$$\varphi(t)\, =\, 1\, -\, \lpa \frac{\sin(\pi\beta)}{\sin(\pi\a)}\, t^{1-\frac{\a}{\beta}}\, +\, \frac{2\sin(\pi(\a-\beta))}{\sin(\pi\a)}\, t\, +\, \frac{\sin(\pi\beta)}{\sin(\pi\a)}\, t^{1+\frac{\a}{\beta}}\rpa\, +\, t^2$$
on $(0,\infty).$ Now since 
$$\frac{\sin(\pi\beta)}{\sin(\pi\a)}\, >\, 0\qquad\mbox{and}\qquad \frac{\sin(\pi(\a-\beta))}{\sin(\pi\a)}\, >\, 0$$
for $1 < \alpha <\beta <2,$ the generalized Descartes' rule of signs - see Theorem 2.2 and Example 2.2 in \cite{Hauk} - implies that $\varphi$ vanishes at most twice on $(0,\infty).$ Let us give some details for completeness: one has
$$\varphi(t)\, =\, \int_0^\infty t^x f(x)\, \mu(dx)$$
where $\mu = \delta_0 + \delta_{1-\frac{\a}{\beta}} + \delta_1 + \delta_{1+\frac{\a}{\beta}} + \delta_2$ is a positive measure, the kernel $(t,x)\mapsto t^x$ is strictly totally positive on $[0,\infty)$ - see Formula (2.1) p.15 in \cite{Karlin}, and the function
$$f\, =\, \Un_{[0,1-\frac{\alpha}{\beta})}\, -\, \lpa \frac{\sin(\pi\beta)}{\sin(\pi\a)}\, \Un_{[1-\frac{\alpha}{\beta},1)}\, +\, \frac{2\sin(\pi(\a-\beta))}{\sin(\pi\a)}\, \Un_{[1,1+\frac{\alpha}{\beta})}\, +\, \frac{\sin(\pi\beta)}{\sin(\pi\a)}\, \Un_{[1+\frac{\alpha}{\beta},2)}\rpa\, +\, \Un_{[2,\infty)}$$
has at most two strict sign changes on $[0,\infty).$ The strict variation-diminishing property - see again Theorem 3.1.(b) p.21 in \cite{Karlin} - implies that $\varphi$ and hence $F_{\a,\b}$ vanishes at most twice on $(0,\infty).$ Moreover, it is clear that $F_{\a,\b}$ is positive in the neighbourhood of zero and $\infty$ and hence must have at least two strict sign changes, since otherwise the function $F_{\alpha}(x) - F_\beta(x)$ would increase and not converge to zero at $\infty$. This shows alltogether that the function $F_{\a,\b}$ vanishes exactly twice on $(0,\infty),$ and completes the proof. 
    
\qed

\begin{Rem}\label{Rem0}

{\em (a) One might ask if the statement of Theorem A extends to higher values of $\a.$ If we first suppose that $\alpha = n$ is an integer and if we set
$$u_n(x)\, =\, E_n(x^n)\, =\, \frac{1}{n}\, \sum_{k=0}^{n-1} e^{xe^{\frac{2\ii k\pi}{n}}}$$
where we have used Formula (3.3.4) in \cite{GKMR} for the second equality, then it can be proved (using $t^4 -4t +3 > 0$ for all $t >1$) that
$$u_3(x)\, =\, \frac{1}{3} \lpa e^x \, +\, 2e^{-x/2} \cos \lpa \frac{\sqrt{3} x}{2}\rpa\rpa \, < \, \cosh(x)\, =\, u_2(x)$$
and (separating the case $x\ge \pi$ and $x\le\pi$ and using $t^3 -3t +2 > 0$ for all $t >1$) that
$$u_4(x)\, =\, \frac{1}{2} \lpa \cosh(x) + \cos(x)\rpa \, <\, u_3(x)$$
for all $x > 0.$ However, using the Wolfram package we found that 
$$u_{25}(30)\, = \, 5.4860...\times 10^{11}\, < \, 6.3108...\times 10^{11} \, = \,u_{26}(30),$$ 
so that the mapping $\alpha \mapsto E_\alpha(x^\alpha)$ has increasing points on $(0,\infty)$ for some $x > 0.$ We believe, and simulations seem to confirm, that the mapping $\a\mapsto\, E_\a(x^\a)$ still decreases on $(2,4)$. Observe that in that case one has the stochastic representation 
\begin{equation}
\label{Ea4}
E_\a(x^\a)\, =\, \frac{1}{2} \lpa E_{\a/2}(x^{\a/2}) \, +\, E_{\a/2}(-x^{\a/2})\rpa\, =\, \frac{1}{2} \lpa \esp\lcr e^{-x \X_{\a/2}}\rcr \, +\, \esp\lcr e^{\ii x \Y_{\a/2}}\rcr\rpa
\end{equation}
where $\X_{\a/2}$ is the centered random variable used in the above proof and $\Y_{\a/2}$ is a symmetric random variable with density
$$\frac{\sin(\pi\alpha/4)\, \vert x\vert^{\alpha/2-1}}{\pi(\vert x\vert^{\alpha}+2\cos(\pi\alpha/4)\vert x\vert^{\alpha/2}+1)}$$ 
on $\rl.$ Above, we have used again Formula (3.3.4) in \cite{GKMR} for the first equality and for the second equality we have applied the Fourier inversion formula to the identity
$$\int_\rl E_{\a/2}(-\vert x\vert^{\a/2})\, e^{\ii ux}\, dx\, = \, 2\,\Re\lpa \frac{(\ii u)^{\a/2-1}}{1+ (\ii u)^{\a/2}}\rpa$$
which is a direct consequence of Formula (3.7.7) in \cite{GKMR}. In view of the latter representation of $E_\a(x^\a)$ as the half-sum of a moment generating function and a characteristic function, it seems that other tools are necessary in order to study the monotonicity of $\a\mapsto E_\a(x^\a)$ for $\a\in (2,4).$ 

\medskip

(b) The mapping $\a\mapsto E_{\a,\b}(x^\a)$ may not be monotonic for $\b < 1.$ For instance, the function
$$E_{\a,\a} (x^\a)\, =\, \a\Ea'(x^\a)\, =\, x^{1-\a}\lpa\frac{e^x}{\alpha}\,+\,\int_0^\infty e^{-xt}\, tf_\a(t)\, dt\rpa$$
is such that 
$$E_{\a,\a} (0)\, -\, E_{\b,\b} (0)\, =\,\frac{1}{\Ga(\a)}\, -\, \frac{1}{\Ga(\b)}\, <\, 0\qquad\mbox{and}\qquad E_{\a,\a} (x^\a)\, -\, E_{\b,\b} (x^\b)\, \to\,\infty\quad\mbox{as $x\to\infty$}$$
if $0<\a< \b < 1.$ If we fix $\b \in (0,1),$ then a consequence of Formula (1.1.16) in \cite{Popov} is the integral representation
$$E_{\a,\b} (x^\a)\, =\, x^{1-\b}\lpa\frac{e^x}{\alpha}\,+\,\int_0^\infty t^{\a -\b} \, e^{-xt}\, \lpa \frac{\sin(\pi\b)\, t^{\alpha}\, +\, \sin(\pi(\a-\b))}{\pi(t^{2\alpha}-2\cos(\pi\alpha)t^\alpha+1)}\rpa dt\rpa$$ for all $\a \in (0,1),$ 
which is however not informative enough at first sight.

\medskip

(c) A combination of \eqref{M-L Repr} and the main result of \cite{Alzer} give also the rigid bounds
\begin{equation}
\label{rigid}
0\, <\, e^x - \b E_\b(x^\b)\, <\, e^x - \a E_\a(x^\a)\, <\, 1-\a
\end{equation}
uniformly on $x > 0.$ Observe that $e^x- \b E_\b(x^\b)\to 0$ as $\b\to 1,$ whereas $e^x- \a E_\a(x^\a)$ converges to a less trivial limit as $\a\to 0.$ More precisely, we can rewrite \eqref{M-L Repr} as
$$e^x- \a E_\a(x^\a)\, =\, (1-\a) \esp\lcr e^{-x \U_\a}\rcr$$
where 
$$\U_\a \, \elaw\, \lpa\frac{\Z_{1-\a}}{\Z_{1-\a}}\rpa^{\frac{1-\a}{\a}}$$
is a positive power transform of the quotient of two independent copies of the positive stable random variable $\Z_{1-\a}$ with Laplace transform $\esp[e^{-t \Z_{1-\a}}] = e^{-t^{1-\a}}$ for all $t\ge 0$ - see Proposition 7 in \cite{Frechet} with $\a = 1$ therein. On the other hand, it follows from the proof of Theorem A pp. 11-12 in \cite{Frechet} that
$$\U_\a \, \claw\, \frac{{\bf S}}{{\bf S}}\qquad \mbox{as $\a\to 0$}$$
where the quotient is again independent and ${\bf S}$ is the log-stable random variable with Mellin transform $\esp[{\bf S}^s] = s^s$ for all $s > 0$ - see the introduction of \cite{Frechet} pp.2-3 for more detail on this random variable. This shows alltogether that
$$e^x- \a E_\a(x^\a)\, \uparrow \, \esp\lcr e^{-x\frac{{\bf S}}{{\bf S}}}\rcr$$
for all $x > 0$ as $\a \to 0.$ We will come back to the boundedness properties of the function $e^x - \a\Ea(x^\a)$ on $(0,\infty)$ for all $\a\in(0,4)$ in Proposition \ref{Bound} below.}

\end{Rem}

We end this section with an application of Theorem A on certain fractional differential equations. We will prove the two following comparison theorems, which seem to have passed unnoticed in the literature on fractional calculus. Our first result deals with an integral equation hinging upon the standard Riemann-Liouville integral
$$I^\a_{0+}f(x)\, =\, \frac{1}{\Ga(\a)}\int_0^x (x-y)^{\a-1} f(y)\, dy$$
for $\a >0,$ which acts on a continuous function $f :\rl^+\to\rl.$ See e.g. Section 1.1 in \cite{AMP} for more details on this fractional integral.

\begin{Coro}
\label{FDE1}
Let $g_i:\rl^+\to\rl$ be a continuous function and $f_i$ be the unique solution to 
$$f_i(x)\, =\, g_i(x) \, +\, \lbd_i^{\a_i} I^{\a_i}_{0+}f(x)$$
with $\lbd_i > 0$ and $\a_i \in (0,1)$ for $i=1,2.$ Suppose $\a_2 \ge \a_1, \lbd_1 \ge \lbd_2$ and $g_1(x)\ge g_2(x)$ for all $x\in\rl^+.$ Then, one has 
$$f_1(x)\,\ge\, f_2(x)$$ 
for all $x\in\rl^+\! .$ 
\end{Coro}

\begin{proof} Applying Formula (7.1.24) in \cite{GKMR}, we can express the two solutions as
$$f_i(x)\, =\, g_i(x)\, +\, \lbd_i^{\a_i} \int_0^x g_i(x-y)\, y^{\a_i-1}\, E_{\a_i,\a_i}\!( (\lbd_i y)^{\a_i})\, dy$$
for $i=1,2$ and all $x\ge 0.$ This implies
$$f_1(x)-f_2(x)\, \ge\, \int_0^x g_1(x-y)\lpa \lbd_1^{\a_1} y^{\a_1-1}\, E_{\a_1,\a_1}\!( (\lbd_1 y)^{\a_1})\, -\, \lbd_2^{\a_2} y^{\a_2-1}\, E_{\a_2,\a_2}\!( (\lbd_2 y)^{\a_2})\rpa dy.$$
On the other hand, it follows from the proof of Theorem A that for all $0 < \a <\b <1$ and $x\ge 0,$ one has
$$E_{\alpha}(x^{\alpha})\,-\,E_{\beta}(x^{\beta})\, = \,\lpa \frac{1}{\alpha}-\frac{1}{\beta}\rpa\lpa e^x\, -\, \esp \lcr e^{-x\X_{\a,\b}}\rcr\rpa$$
where $\X_{\a,\b}$ is a positive random variable with density function
$$t\,\mapsto\,\frac{\a\b (f_\a(t) - f_\b(t))}{\b-\a}\cdot$$
In particular, the function 
\begin{eqnarray*}
E_{\a_1}\!((\lbd_1 y)^{\a_1}) - E_{\a_2}\!((\lbd_2 y)^{\a_2}) & = & E_{\a_1}\!((\lbd_1 y)^{\a_1}) - E_{\a_1}\!( (\lbd_2 y)^{\a_1})\, + \, E_{\a_1}\!((\lbd_2 y)^{\a_1}) - E_{\a_2}\!((\lbd_2 y)^{\a_2})\\
& = & \sum_{n\ge 0} \frac{(\lbd_1^{n\a_1} -\lbd_2^{n\a_1}) y^{n\a_1}}{\Ga(1+n\a_1)}\, +\, \lpa \frac{1}{\alpha_1}-\frac{1}{\a_2}\rpa\lpa e^{\lbd_2 y}- \esp \lcr e^{-\lbd_2 y\X_{\a_1,\a_2}}\rcr\rpa
\end{eqnarray*}
increases on $(0,\infty).$ Taking the derivative shows that
\begin{equation}
\label{Derive}
\lbd_1^{\a_1} y^{\a_1-1}\, E_{\a_1,\a_1}\!( (\lbd_1 y)^{\a_1})\, \ge\, \lbd_2^{\a_2} y^{\a_2-1}\, E_{\a_2,\a_2}\!( (\lbd_2 y)^{\a_2})
\end{equation}
for all $y\ge 0,$ and concludes the proof.
\end{proof}

Our second result deals with the standard Riemann-Liouville derivative
$$D^\a_{0+}f(x)\, =\, \frac{1}{\Ga(1-\a)}\frac{d}{dx}\lpa\int_0^x (x-y)^{-\a} f(y)\, dy\rpa $$
for $\a \in(0,1),$ which acts on an absolutely continuous function $f :\rl^+\to\rl.$ See e.g. Section 1.7 in \cite{AMP} for more details on this fractional derivative.

\begin{Coro}
\label{FDE2}
Let $g_i:\rl^+\to\rl$ be a continuous function and $f_i$ be the unique solution to 
$$\lacc\begin{array}{l}
D^{\a_i}_{0+} f_i(x)\, -\, \lbd f_i(x) \; =\; g_i(x)\\
I^{1-\a_i}_{0+} f_i(0)\, =\, \mu_i\end{array}\right.$$
with $\lbd\ge 1, \mu_i\ge 0$ and $\a_i \in (0,1)$ for $i=1,2.$ Suppose $\a_2 \ge \a_1, \mu_1\ge \mu_2$ and $g_1(x)\ge g_2(x)$ for all $x\in\rl^+.$ Then, one has 
$$f_1(x)\,\ge\, f_2(x)$$ 
for all $x\in\rl^+\! .$ 
\end{Coro}

\begin{proof}
Applying Formula (2.96) in \cite{AMP}, we can express the two solutions as
$$f_i(x)\, =\, \mu_i\, x^{\a_i-1} E_{\a_i,\a_i}\; (\lbd x^{\a_i})\, +\, \int_0^x g_i(x-y)\, y^{\a_i-1}\, E_{\a_i,\a_i}(\lbd y^{\a_i})\, dy$$
for $i = 1,2.$ Defining $\lbd_1, \lbd_2\ge 1$ by $\lbd_1^{\a_1} = \lbd_2^{\a_2} = \lbd,$ we have $\lbd_1 \ge \lbd_2$ since $\lbd\ge 1$ and $\a_2 \ge \a_1,$ and deduce from \eqref{Derive} that 
$$t^{\a_1-1}\,E_{\a_1,\a_1}( \lbd t^{\a_1})\, \ge\, t^{\a_2-1}\, E_{\a_2,\a_2}( \lbd t^{\a_2})$$
for all $t\ge 0.$ This completes the proof.
\end{proof}

\begin{Rem}
\label{RemFDE}
{\em Applying Proposition 2.61 in \cite{AMP}, we see that the same comparison theorem holds for the fractional Cauchy problem
$$\lacc\begin{array}{l}
^C\! D^{\a}_{0+} f(x)\, -\, \lbd f(x) \; =\; g(x)\\
f(0)\, =\, \mu\end{array}\right.$$
with $\mu\ge 0$ and $\lbd\ge 1,$ where $^C\! D^{\a}_{0+}$ is the Dzhrbashyan-Caputo fractional derivative. Indeed, the above proof shows that we also have $E_{\a_1}( \lbd t^{\a_1})\ge E_{\a_2}( \lbd t^{\a_2})$ for all $\lbd\ge 1$ and $t\ge 0$ as soon as $0 < \a_1\le\a_2 < 1.$ Observe that the same argument shows that the comparison theorem also holds for the integral equation
$$f(x)\, =\, g(x) \, +\, \lbd I^{\a}_{0+}f(x)$$
with $\lbd\ge 1$ and $\a \in (0,1).$}
\end{Rem}

\section{Proof of Theorem B} 

We start with the standard formula which is given in e.g. in Exercise 3.9.5. of \cite{GKMR}:
\begin{equation}
\label{ML Rep}
E_\a(-x^\a)\, =\,\int_0^\infty e^{-xt}g_\a(t)\,dt, \qquad x\ge 0,
    \end{equation}
with
$$g_\a(t)\,=\, \frac{\sin(\pi\alpha)\, t^{\alpha-1}}{\pi(t^{2\alpha}+2\cos(\pi\alpha)t^\alpha+1)}$$ 
for all $t > 0.$ Beware in passing that this formula is valid for $\a\in (0,1)$ only, since the function $x\mapsto E_\a(-x^\a)$ vanishes at least once on $(0,\infty)$ for every $\a > 1$ - see e.g. \cite{Poly, Popov}. Setting $\V_\a$ for the random variable with density $g_\a$ and $G_{\a,\b}(x) = E_\a(-x^\a) -E_\b(-x^\b),$ an integration by parts gives
$$G_{\a,\b}(x)\, =\,  x\int_0^\infty e^{-xt} \lpa \pb [ \V_\b > t] - \pb [ \V_\a > t]\rpa dt.$$ On the other hand, the function $t\mapsto \psi_{\a,\b} (t) = \pb [ \V_\b > t] - \pb [ \V_\a > t]$ has derivative
$$g_\a(t)\,-\, g_\b(t)\,=\, \frac{\sin(\pi\alpha)\, t^{\alpha}(1+ t^{2\beta}) \, - \, \sin(\pi\beta)\, t^{\beta}(1+ t^{2\a}) \, +\, 2\sin(\pi(\a-\beta))\, t^{\a+\beta}}{\pi\,t(t^{2\a} + 2\cos(\pi\a) t^\a +1)(t^{2\beta} + 2\cos(\pi\beta) t^\beta +1)}$$ 
and the fact that $\sin(\pi(\a-\b)) < 0$ shows as in the proof of Theorem A that it has at most two strict sign changes, starting positive. Hence, the function $\psi_{\a,\b}$ which vanishes at zero and infinity, has exactly one strict sign change on $(0,\infty),$ starting positive. Putting this together with Theorem 3.1.(b) p.21 in \cite{Karlin} implies that the function $G_{\a,\b}(x)$ vanishes at most once on $(0, \infty).$ Finally, since
$$G_{\a,\b}(x)\,\sim\, \frac{1}{\Gamma(1-\a)\, x^\a}\, > \, 0 \quad \mbox{as $x\to\infty$}\qquad\mbox{and}\qquad
G_{\a,\b}(x)\,\sim\, -\frac{x^\a}{\Gamma(1+\a)}\, < \, 0 \quad \mbox{as $x\to 0,$}$$ 
we have shown that there exists $x_{\a,\b} > 0$ with the required properties. 

The lower bound on $x_{\a,\b}$ is a consequence of Theorem 4 in \cite{Frechet}, which gives
$$G_{\a,\b}(x)\;\le\; \frac{\Ga(1-\b)\Ga(1+\a)x^\b -x^\a}{(\Ga(1+\a) + x^\a)(1+\Ga(1-\b)x^\b)}\; <\; 0$$ 
as soon as $x < (\Ga(1-\b)\Ga(1+\a))^{-\frac{1}{\b-\a}}.$ This result also implies
$$G_{\a,\b}(x)\;\ge\; \frac{x^\b - \Ga(1-\a)\Ga(1+\b) x^\a}{(\Ga(1+\b) + x^\b)(1+\Ga(1-\a)x^\a)}\; > \; 0$$ 
for $x > (\Ga(1-\a)\Ga(1+\b))^{\frac{1}{\b-\a}},$ which yields the upper bound $x_{\a,\b}  < \lpa\Ga(1-\a)\Ga(1+\b)\rpa^{\frac{1}{\b-\a}}.$
The latter is however weaker than the one we claim, because $\Ga(1+\a)\Ga(1-\a) = 1/\sin_c(\a\pi) > 1.$ In order to get the required upper bound, we will first use a fractional differential equation argument. Specifically, it follows from Exercise 3.9.4. in \cite{GKMR} and a change of variable that 
$$G_{\a,\b}(x)\; =\; \frac{x^\b}{\Ga(1+\b)}\int_0^1 \b(1-t)^{\b-1} E_\b (-(tx)^\b)\, dt \, -\, \frac{x^\a}{\Ga(1+\a)}\int_0^1 \a(1-t)^{\a-1} E_\a (-(tx)^\a)\, dt.$$  
If we now set
$$x_\ast \, =\, \lpa\frac{\Ga(1+\b)}{\Ga(1+\a)}\rpa^{\frac{1}{\b-\a}}$$
and if we suppose $x_\ast \le x_{\a,\b},$ the definitions of $x_\ast$ and $x_{\a,\b}$ imply
\begin{eqnarray*}
0 \; \ge\; G_{\a,\b}(x_\ast)  & = & \frac{x_\ast^\b}{\Ga(1+\b)} \lpa \int_0^1 \b(1-t)^{\b-1} E_\b (-(tx_\ast)^\b)\, dt \, -\, \int_0^1 \a(1-t)^{\a-1} E_\a (-(tx_\ast)^\a)\, dt\rpa \\
& > & \frac{x_\ast^\b}{\Ga(1+\b)}\int_0^1 \lpa \b(1-t)^{\b-1} -\a(1-t)^{\a-1}\rpa E_\a (-(tx_\ast)^\a)\, dt\\
& = & \frac{x_\ast^\b}{\Ga(1+\b)} \lpa \esp \lcr \Ea\lpa -\lpa x_\ast\B_{1,\b}\rpa^\a)\rpa\rcr \, - \, \esp \lcr \Ea\lpa -\lpa x_\ast\B_{1,\a}\rpa^\a)\rpa\rcr\rpa\; > \; 0
\end{eqnarray*}
where in the second equality we have introduced the standard beta random variable and the second inequality follows from the easily established fact that $\B_{1,\b}\prec_{st} \B_{1,\a}$ where $\prec_{st}$ stands for the usual stochastic order, the decreasing character of the function $t\mapsto \Ea \lpa-\lpa tx_\ast\rpa^\a\rpa$ on $(0,1),$ and Theorem 1.A.3.(a) in \cite{Shaked}. This contradiction implies $x_{\a,\b} < x_\ast$ as required. 

We finally show that $x_{\a,\b} < 1,$ which amounts to $G_{\a,\b}(1) > 0$ with the above notation, and for this we will use yet another argument based on the useful property that $\V_\a$ and $\V_\b$ are self-reciprocal random variables, that is
$$\V_\a\,\elaw\,\frac{1}{\V_\a}\qquad\mbox{and}\qquad \V_\b \, \elaw\, \frac{1}{\V_\b}\cdot$$
The latter implies $\psi_{\a,\b} (u) +\psi_{\a,\b}(1/u) = 0$ for all $u\in (0,\infty)$ with the above notation, so that $\psi_{\a,\b}(1) = 0$ and $\psi_{\a,\b}$ is positive on $(0,1).$ Therefore, 
$$G_{\a,\b}(1)\, = \, \int_0^\infty e^{-u} \psi_{\a,\b} (u)\, du \, =\, \int_0^1  \lpa e^{-u} \, -\, u^{-2} e^{-1/u}\rpa\, \psi_{\a,\b} (u)\, du\, > \,0$$    
where the last strict inequality comes from the positivity of the integrand (indeed, for all $u\in (0,1)$ one has $\rho (u) = u^2e^{1/u -u} > 1$ since $\rho(1) = 1$ and $(\log\rho)'(u) < 0)$. This completes the proof.

\qed

\begin{Rem} 

\label{Rem1}{\em (a) One might wonder whether there is a more direct proof of the decreasing character on $(0,1)$ of the mapping 
$$\a\,\mapsto \,\Ea(-1)\, =\, \sum_{n\ge 0} \, \frac{(-1)^n}{\Ga(1+\a n)}$$ 
This is equivalent to the positivity of the series
$$\sum_{n\ge 1}\, (-1)^n \lpa\frac{\psi(1+\a n)}{\Ga(\a n)}\rpa$$  
for all $\a\in (0,1),$ but we could not find any simple argument for that.

\medskip

(b) The bound $x_{\a,\b} < 1$ shows that for all $x\ge 1,$ the mapping $\a\mapsto\Ea(-x^\a)$ decreases on $(0,1),$ and this yields another proof of the decreasing character on $(1,2)$ of the mapping $\a\mapsto \Ea(x^\a)$ for all $x > 0.$ Indeed, it implies together with Formula 3.3.4. in \cite{GKMR} and the case $\a\in (0,1)$ that the mapping decreases for all $x\ge 1.$ On the other hand, the log-concave function
$$\a\mapsto \frac{e^{-ts\a}}{\Ga(1+t\a)}$$ 
has logarithmic derivative $-t(s +\psi(1+t)) < 0$ at $\a = 1$ for all $s > 0$ and $t\ge 1,$ so that the mapping
$$\a\mapsto \frac{x^{n\a}}{\Ga(1+n\a)}$$ 
decreases on $(1,\infty)$ for all $x < 1,$ which readily implies that $\a\mapsto\Ea(x^\a)$ decreases on $(1,\infty)$ for all $x < 1$ as well, and completes the argument. All in all, this second proof also relies on convex ordering, albeit in a less direct way than the one we used for the proof of Theorem A.

\medskip

(c) It follows from Theorem 4 in \cite{Frechet} that $\Ea(-x^\a)\to 1/2$ as $\a\to 0$ for all $x > 0,$ so that the statement of Theorem B remains true for $\a = 0$ and $\b\in (0,1]$ with $x_{0,\b}$ defined as the unique solution to $E_\b (-x^\b) = 1/2$ on $(0,\infty),$ in other words $x_{0,\b} = (-\log_\b(1/2))^{1/\b}$ where $\log_\b$ is the generalized logarithm introduced in our previous paper \cite{FS2}. Observe that Theorem B remains also true for $\b = 1$ and $\a\in (0,1),$ as a consequence of the log-convexity of the function $x\mapsto e^x \Ea(-x^\a)$ which is decreasing-then-increasing, starts from 1 and tends to $\infty$ as $x\to\infty.$ Observe finally that our bounds yield $x_{\a,\b}\to e^{-\g}$ as $0<\a <\b\downarrow 0.$}
\end{Rem}

The following consequence of Theorem B is a similar sign-change property for the functions $x\mapsto \Ea(-\lbd x^\a)$ with $\lbd > 0.$ Apart from being the fundamental solution to a contracting Abelian integral equation - see again Exercise 3.9.4 in \cite{GKMR}, these functions are important when $\lbd$ is the eigenvalue of a certain diffusion operator - see \cite{LMS} for a pioneering work on the subject and in particular Theorem 3.2 therein.
  
\begin{Coro}
\label{Iden}
For every $0 <\a <\b <1$ and $\lbd > 0,$ there exists a unique $x_{\a,\b,\lbd} > 0 $ such that the function 
$$x\,\mapsto\,\Ea(-\lbd x^\a) - E_\b(-\lbd x^\b)$$ 
is negative for $x\in (0,x_{\a,\b,\lbd})$ and positive for $x\in (x_{\a,\b,\lbd},\infty).$
\end{Coro}

\begin{proof}
The existence of $x_{\a,\b,\lbd}$ is established similarly as for Theorem B via the formula
$$G_{\a,\b,\lbd} (x)\, =\, \Ea(-\lbd x^\a) - E_\b(-\lbd x^\b)\, =\, \lbd \int_0^\infty e^{-xt} \,g_{\a,\b,\lbd} (t)\, dt,$$ 
where
$$g_{\a,\b,\lbd} (t) \, =\, \frac{\sin(\pi\alpha)\, t^{\alpha}(\lbd^2 + t^{2\beta}) \, - \, \sin(\pi\beta)\, t^{\beta}(\lbd^2 + t^{2\a}) \, +\, 2\lbd^2 \sin(\pi(\a-\beta))\, t^{\a+\beta}}{\pi\,t\,(t^{2\a} + 2\lbd\cos(\pi\a) t^\a +\lbd^2)(t^{2\beta} + 2\lbd\cos(\pi\beta) t^\beta +\lbd^2)}$$
has again at most two strict sign changes on $(0,\infty)$. Hence, the function $G_{\a,\b,\lbd}(x)$ vanishes at most twice on $(0,\infty)$ and since it is negative in the neighbourhood of zero and positive in the neighbourhood of $\infty,$ it vanishes once as required. 

\end{proof}

\begin{Rem} 
\label{Niko}{\em If $\lbd < 1,$ choosing $x =\lbd^{-1/\a}$ we have
$$\Ea(-\lbd x^\a)\, = \, \Ea(-1) \, > \, E_\b(-1)\, >\, E_\b ( -\lbd^{1-\b/\a}) \, = \, E_\b(-\lbd x^\b)$$ 
where the first inequality comes from Remark \ref{Rem1} (a) and the second inequality from $\lbd^{1-\b/\a} > 1$ and the decreasing character of $u \mapsto E_\b(-u)$ on $(0,\infty).$ This implies the upper bound 
$$x_{\a,\b,\lbd}\,\le\, \lbd^{-1/\a}.$$ 
If $\lbd > 1$ and if we further suppose $\b\le\b_0 = 0.46163..$ where $\Ga(1+\b_0)$ is the minimum of the Gamma function on $(0,\infty),$ a direct consequence of Theorem C is that 
$$\Ea(-\lbd) \, =\, \Ea\lpa \Ga(1+\a)\lpa\frac{-\lbd}{\Ga(1+\a)}\rpa\rpa\, > \,  E_\b\lpa \Ga(1+\b)\lpa\frac{-\lbd}{\Ga(1+\a)}\rpa\rpa\, >\, E_\b(-\lbd)$$ 
for all $\a < \b,$ which implies $x_{\a,\b,\lbd} \le 1.$ We believe that the upper bound
$$x_{\a,\b,\lbd}\,\le\, 1$$
holds for all $\lbd > 1$ and $0 < \a<\b < 1.$ This would show the interesting property that the mapping $\a\mapsto E_\a ( -\lbd x^\a)$ decreases on $(0,1)$ for all $\lbd, x\ge 1.$ See Lemma 6 in \cite{AZ} for a related result, based on a contour argument.}
\end{Rem}

The next proposition is related to Theorem 3.6 in \cite{AS}. We recall that $x_{\a,\b}$ is the unique zero of the function $\Ea(-x^\a)\, - \, E_\b(-x^\b)$ on $(0,\infty).$ 

\begin{Propo}
\label{Eaa1}
For every $0 <\a <\b <1,$ there exists $y_{\a,\b} < x_{\a,\b}$ and $z_{\a,\b} > x_{\a,\b}$ such that the function 
$$x\,\mapsto\, x^{\b-1}E_{\b,\b}(-x^\b)\, - \,x^{\a-1}E_{\a,\a}(-x^\a)$$ 
is negative for $x\in (0,y_{\a,\b})\cup (z_{\a,\b},\infty)$ and positive for $x\in (y_{\a,\b},z_{\a,\b}).$ 
\end{Propo}

\begin{proof} With the notations of the proof of Theorem B, one has
$$x^{\b-1}E_{\b,\b}(-x^\b)\, - \,x^{\a-1}E_{\a,\a}(-x^\a)\, =\, G_{\a,\b}'(x)\, =\,  \int_0^\infty e^{-xt}\, t\lpa g_\b(t) - g_\a(t)\rpa dt$$  
and since the function $t\mapsto t \lpa g_\b(t) - g_\a(t)\rpa$ has two strict sign-changes on $(0,\infty),$ the function $G_{\a,\b}'(x)$ vanishes at most twice. On the other hand, one has
$$G_{\a,\b}'(x)\, \sim\, - \frac{x^{\a-1}}{\Ga(\a)}\, <\, 0\quad\mbox{as $x\to 0$} \qquad\mbox{and}\qquad G_{\a,\b}'(x)\, \sim\, \frac{x^{-\a-1}}{\Ga(-\a)}\, <\, 0\quad\mbox{as $x\to \infty,$}$$ 
and since the function $G'_{\a,b}$ must vanish at least once on $(0,\infty)$ as the difference of two densities, this implies that it vanishes exactly twice as required. Setting $y_{\a,\b} < z_{\a,\b}$ for the two roots, it is clear from the fact that 
$$G_{\a,\b} (0)\, = \, G_{\a,\b}(x_{\a,\b})\, =\, \lim_{x\to\infty} G_{\a,\b} (x)\, = \,0,$$ 
combined with Rolle's theorem, that one must have $y_{\a,\b} < x_{\a,\b}< z_{\a,\b}.$

\end{proof}

\begin{Rem}\label{RemEa}
{\em (a) This result has a direct application on the Mittag-Leffler random variable ${\hat \M}_\a$ with parameter $\a\in (0,1),$ which had been introduced in \cite{Pill} as the positive random variable with Laplace transform
$$\esp\lcr e^{-\lbd {\hat \M}_\a}\rcr\, =\, \frac{1}{1+\lbd^\a}$$
for $\lbd\ge 0.$ Indeed, it follows from Theorem 2.1 in \cite{Pill} that ${\hat \M}_\a$ has density $x^{\a-1}E_{\a,\a}(-x^\a)$ on $(0,\infty),$ and from Theorem 2.2 in \cite{Pill} and logarithmic differentiation at $\rho = 0$ that $\log{\hat \M}_\a$ is integrable with $\esp[\log{\hat \M}_\a] = -\gamma$ for all $\a\in (0,1),$ where $\gamma$ is the Euler-Mascheroni constant. Therefore, since the density of $\log{\hat \M}_\a$ on $\rl$ is $e^{\a x}E_{\a,\a}(-e^{\a x}),$ we deduce from Proposition \ref{Eaa1} and Theorem 3.A.44. in \cite{Shaked} that
$$\log {\hat \M}_\b\, \prec_{cx} \, \log {\hat \M}_\a$$
for all $0 <\a<\b < 1.$ Observe that $\log {\hat \M}_\b$ converges in distribution towards $\log \L$ as $\b \to 1,$ where $\L$ is the unit exponential random variable, whereas $\log {\hat \M}_\a$ does not converge in distribution as $\a \to 0.$ The distribution function of $\log{\hat \M}_\a$ is $1- \Ea(-e^{\a x}),$ and it is worth comparing this random variable with the fractional Gumbel random variable $\G_{\a,\lbd}$ recently introduced in \cite{BSV} and having distribution function $1- \cL_\a(-e^{\lbd x}),$ where 
$$\cL_\a(x)\, =\, \sum_{n\ge 0} \frac{x^n}{(n!)^\a}$$ 
is the classical Le Roy function - see Theorem 1.3 therein. It follows indeed from Proposition 10 in \cite{BS} that 
$$\G_{\b,\lbd}\, \prec_{st} \, \G_{\a,\lbd}$$
for all $0 <\a<\b < 1$ and $\lbd\ge 0,$ where $\prec_{st}$ stands for the usual stochastic order. 

\medskip

(b) Using the bounds on $E_{\a,\a}(x)$ given in Remark \ref{Rem4} (a) below, a lower bound on $y_{\a,\b}$ resp. an upper bound on $z_{\a,\b}$ can be given as the smallest positive root resp. largest positive root to some polynomial equation with fractional powers. These two roots are however non explicit. Simulations also suggest that the mapping $\a\mapsto E_{\a,\a}(-1)$ increases on $(0,1),$ which would imply the interesting bound 
$$x_{\a,\b} \, < \,1 \,<\, z_{\a,\b}.$$ 
Similarly as in Remark \ref{Rem1} (a), this amounts to the positivity of the alternate series
$$\sum_{n\ge 1} (-1)^n \lpa\frac{n \psi(\a n)}{\Ga(\a n)}\rpa\, =\, \frac{1}{\a} \, E_{\a,\a}(-1)\, +\, \sum_{n\ge 1} (-1)^n \lpa\frac{n \psi(1+\a n)}{\Ga(\a n)}\rpa$$  
for all $\a\in (0,1),$ but we were unfortunately unable to prove this.}
\end{Rem}

Our next proposition gives a explanation why the main result of \cite{Alzer} cannot hold in general. See in particular the comment after the statement of the Theorem therein.
 
\begin{Propo}
\label{SC2}
For all $1 <\a <\b <2,$ there exists $x_{\a,\b}^\ast > x_{\a/2,\b/2}$ such that the function 
$$x\,\mapsto\,\b E_\b(x^\b)\, - \,\a\Ea(x^\a)$$ 
is positive for $x\in (0,x_{\a,\b}^\ast)$ and negative for $x\in (x_{\a,\b}^\ast,\infty).$ 
\end{Propo}

\begin{proof} Setting $H_{\a,\b}(x) = \b E_\b(x^\b)\, - \,\a\Ea(x^\a)$ for all $x\ge 0,$ we have the formula
$$H_{\a,\b}(x) \, =\, \int_0^\infty e^{-xt} \lpa\frac{\a\sin(\pi\a) t^{\a -1}}{\pi (t^{2\a} - 2\cos(\pi\a) t^\a +1)}\, -\, \frac{\b\sin(\pi\b) t^{\b -1}}{\pi (t^{2\b} - 2\cos(\pi\b) t^\b +1)}\rpa dt$$
as a direct consequence of \eqref{M-L Rep}. In particular, an easy asymptotic analysis shows that
$$H_{\a,\b}(x) \, \sim\, -\frac{x^{-\a}}{\Ga(-\a)}\, <\, 0\qquad\mbox{as $x\to\infty$}$$
and since $H_{\a,\b}(x) = \b-\a > 0,$ this implies that the function $H_{\a,\b}$ has at least one strict sign change on $(0,\infty).$ Another consequence of the differentiation of \eqref{M-L Rep} at $x=0$ is that the function 
$$h_\a(t)\,=\,\frac{(-\a\sin(\pi\a)) t^\a}{\pi (t^{2\a} - 2\cos(\pi\a) t^\a +1)}$$
is a density function on $(0,\infty)$ for all $\a\in (1,2).$ If we set ${\tilde \X_\a}$ for the corresponding random variable, we have
$$H_{\a,\b}'(x) \, =\, \esp\lcr e^{-x{\tilde \X_\a}}\rcr\, -\, \esp\lcr e^{-x{\tilde \X_\b}}\rcr\, =\, x\int_0^\infty e^{-xt} (\pb [ {\tilde \X_\b} > t] - \pb [ {\tilde \X_\a} > t])\,dt$$
for all $1 <\a <\b <2.$ On the other hand, the function
$$h_\a(t)\,-\,h_\b(t)\, =\, \frac{(-\a\sin(\pi\alpha))\, t^{\alpha}(1+ t^{2\beta}) \, + \, \b\sin(\pi\beta)\, t^{\beta}(1+ t^{2\a}) \, +\, 2\psi(\a,\b)\, t^{\a+\beta}}{(t^{2\a} - 2\cos(\pi\a) t^\a +1)(t^{2\beta} - 2\cos(\pi\beta) t^\beta +1)}$$
with $\psi(\a,\b) = \a\sin(\pi\a)\cos(\pi\b)-\b\sin(\pi\b)\cos(\pi\a) < 0$ has at most two strict sign changes on $(0,\infty),$ starting positive, for the same reason as in the proof of Theorem B. The negative character of $\psi(\a,\b)$ is clear for $1 <\a\le 3/2$ and $\a <\b <2$ since 
$$\psi(\a,\b)\, =\, \a\sin(\pi(\a-\b)) + (\a-\b)\sin(\pi\b)\cos(\pi\a),$$ 
whereas for $3/2 < \a <\b < 2$ it comes from $\psi(\a,\b) = \cos(\pi\a)\cos(\pi\b) (\a\tan(\pi\a) -\b\tan(\pi\b))$ and the increasing character of $x\mapsto x\tan(\pi x)$ on $(3/2,2).$ We can then conclude as in the proof of Theorem B to show the existence of $x_{\a,\b}^\ast > 0$ with the required properties. 

Finally, the lower bound $x_{\a,\b}^\ast > x_{\a/2,\b/2}$ comes from the decomposition
$$H_{\a,\b}(x)\, =\, H_{\a/2,\b/2}(x)\, +\, \frac{\b}{2}\, G_{\b/2,\a/2}(x)\, +\, \lpa\frac{\b-\a}{2}\rpa E_{\a/2}(-x^{\a/2})$$
with the notation of the proof of Theorem B, which implies together with the main result of \cite{Alzer} and the positive character of $E_{\a/2}(-x^{\a/2})$ that
$$H_{\a,\b}(x)\, > \, \frac{\b}{2}\, G_{\b/2,\a/2}(x)\,\ge \, 0$$
for all $x \le x_{\a/2,\b/2},$ by Theorem B.

\end{proof}

\begin{Rem}\label{Rem2}

{\em (a) Contrary to Theorem B, we could not get any relevant upper bound on $x^\ast_{\a,\b}.$ Observe that $x^\ast_{\a,\b}\to \infty$ for all $\b\in (1,2)$ as $\a\to 1,$ because of \eqref{M-L Rep}. 

\medskip

(b) Another consequence of \eqref{M-L Rep} is that the statement of Proposition \ref{SC2} remains true for $\beta = 2$ and every $\a\in (1,2),$ where $\W_\a$ is the positive random variable with density 
$$\frac{\a}{\a-1}\lpa\frac{(-\sin(\pi\a)) t^{\a -1}}{\pi (t^{2\a} - 2\cos(\pi\a) t^\a +1)}\rpa$$ 
and $x^\ast_{\a,2}$ is defined as the positive solution - which is again unique by log-convexity as in Remark \ref{Rem1} (c) - to
$$e^x\,\esp[e^{-x\W_\a}] \, =\, \frac{1}{\a -1}\cdot$$}     
\end{Rem}

Our last proposition in this section connects the function $E_\a(x^\a)$ and the leading term in its asymptotic expansion at $\infty$. This result is certainly well-known to practitioners of the Mittag-Leffler function, even though we could not find its precise statement in the literature.

\begin{Propo}
\label{Bound}
The function $x\mapsto \a\Ea(x^\a) - e^x$ is bounded on $(0,\infty)$ if and only if $\a\le 4.$ Besides, it converges to zero as $x\to \infty$ if and only if $\a < 4.$ 
\end{Propo}

\begin{proof} The only if part comes from the asymptotic behaviour
$$\a\Ea(x^\a)\, -\, e^x\; = \; 2e^{x\cos(2\pi/\a)}\cos(x\sin(2\pi/\a)) \, +\, {\rm o}\lpa e^{x\cos(2\pi/\a)}\rpa\qquad\mbox{as $x\to\infty$}$$
for $\a > 4,$ which is given in Formula (3.4.16) of \cite{GKMR}. This formula also gives the if part but we prefer to present an argument compiling our previous discussions, because it leads to optimal bounds. Setting $G_\a(x) = e^x -\a\Ea(x^\a),$ it follows indeed from \eqref{M-L Repr} and \eqref{rigid} that for $\a\in(0,1)$ 
$$G_\a(x) \,\in\, (0, 1-\a)\quad\mbox{for all $x\in (0,\infty)$}\quad\mbox{with}\quad \lim_{x\to 0} G_\a(x) = 1-\a\quad\mbox{and}\quad \lim_{x\to \infty} G_\a(x) = 0,$$  
and from \eqref{M-L Rep} that for $\a\in (0,2)$ one has $G_\a(x) = (1-\a) \esp \lcr e^{-x \W_\a}\rcr  \in (1-\a,0)$
for all $x > 0$ with the above notation for the random variable $\W_\a,$ 
so that $\lim_{x\to 0} G_\a(x) = 1-\a$ and $\lim_{x\to \infty} G_\a(x) = 0.$ Moreover, the decomposition \eqref{Ea4} implies 
$$G_\a(x) \,\in\, (1-\a,\a/2)\quad\mbox{for all $x\in (0,\infty)$}\quad\mbox{with}\quad \lim_{x\to 0} G_\a(x) = 1-\a\quad\mbox{and}\quad \lim_{x\to \infty} G_\a(x) = 0$$  
for $\a\in (2,4),$ the limit at infinity being a consequence of the Riemann-Lebesgue lemma. Finally, the cases $\a=1,2$ are obvious with $G_1(x) = 0$ and $G_2(x) = e^{-x},$ whereas for $\a = 4$ one has
$$G_4(x) \, =\, -\lpa e^{-x} \, +\, 2\cos(x)\rpa,$$
which belongs to $(-3,2)$ on $(0,\infty)$ but does not converge to zero as $x\to\infty.$ 

\end{proof}

\section{Proof of Theorem C} 

We begin with the classical representation $\Ea(x) = \esp[e^{x\M_\a}]$ where $\M_\a$ is the Mittag-Leffler random variable. See (5) in our previous paper \cite{FS2} and also the introduction and the references therein for more material on this random variable. Setting $H_\a(x) = E_\alpha(\Ga(1+\a) x),$ it follows from (11) in \cite{FS2} that $H_\a(x) = \esp[e^{x\T_\a}]$ where $\T_\a$ is defined as the following a.s. convergent infinite independent product of renormalized beta random variables:
$$\T_\a\,=\, \prod_{n=0}^\infty \lpa \frac{n+1}{n+\a}\rpa \B_{1+\frac{n}{\a}, \frac{1}{\a}-1}.$$
We will now show that the mapping $\a\mapsto \T_\a$ is non-increasing for the convex ordering, which will show that the mapping $\a\mapsto H_\a(x)$ is non-increasing for all $x\in\rl$ and hence decreasing for all $x\in\rl^\ast$ since it is non-constant and analytic. To obtain the convex ordering, we use the alternative representation
$$\T_\a\,\elaw\, \prod_{n=0}^\infty \lpa \frac{\Ga(1+(n+1)\a)\Ga((n+1)\a)}{\Ga(1+n\a)\,\Ga((n+2)\a)}\rpa \B_{\a(1+n), 1-\a}^{\a}$$
which is obtained in comparing (3.3) and (3.4) in \cite{Sun} with $\gamma = 1.$ From the stability of the convex ordering by mixtures - see Corollary 3.A.22 in \cite{Shaked}, it is enough to show that the mapping 
$$\a\, \mapsto\, \lpa \frac{\Ga(1+(n+1)\a)\Ga((n+1)\a)}{\Ga(1+n\a)\,\Ga((n+2)\a)}\rpa \B_{\a(1+n), 1-\a}^{\a}$$
decreases on $(0,1)$ for the convex ordering for every $n\ge 0.$ We will show the more general fact that 
the mapping
$$\a\; \mapsto\; \T_\a^{a,b}\; =\; \lpa \frac{\Ga(1+ (b+1)\a) \,\Ga(1+(a+b)\a)}{\Ga(1+b\a)\,\Ga(1+ (a+b+1)\a)}\rpa \B_{\a(1+b), 1-\a}^{\a a}$$
decreases for the convex ordering on $(0,1)$ for all $a,b > 0.$ Indeed, this implies our claim taking $a=1, b=n$ and using
$$\lpa \frac{\Ga(1+(n+1)\a)\Ga((n+1)\a)}{\Ga(1+n\a)\,\Ga((n+2)\a)}\rpa\, =\, \frac{n+2}{n+1}\lpa \frac{\Ga(1+(n+1)\a)\,\Ga(1+(n+1)\a)}{\Ga(1+n\a)\,\Ga(1+(n+2)\a)}\rpa.$$
To do so, we first show that the prefactor giving the right end of the support is a decreasing function in $\a$. Taking the logarithmic derivative and reparametrizing, this amounts to 
$$c\psi(1+c) \,+\, f \psi (1+ f)\, >  \,d \psi (1+d)\, + \, e\psi(1+e)$$ 
for all $0 < c < d$ and $e < f$ with $c+f = d+e.$ Recalling
$$\psi(1+z)\; = \; \psi(1) \, +\, \sum_{n\ge 1} \frac{z}{n(n+z)}$$
for all $z > 0,$ this is a consequence of the strict convexity of $x\mapsto x^2/(n+x)$ on $(0,\infty)$ for all $n\ge 1.$ Setting $t_\a(x)$ for the density of $\T_\a^{a,b}$ over the interior of its support $(0,m_\a),$ we deduce that $t_\a(x) > t_\b (x) = 0$ if $x\in ]m_\b,m_\a[,$ whereas $t_\a(x) = t_\b (x) = 0$ if $x > m_\a$ and $t_\b(x)\to \infty > t_\a(m_\b)$ as $x\uparrow m_\b,$ for all $0 < \a < \b < 1.$ Moreover the densities $t_\a(x)$ and $t_\b(x)$ cross at least once on $]0,m_\b[$ since otherwise we would have 
$$\T_\b^{a,b}\, \prec_{st}\, \T_\a^{a,b}$$ 
by Theorem 1.A.12 in \cite{Shaked}, which is impossible by the equality of expectations. Using again Theorem 3.A.44 in \cite{Shaked}, we are reduced to prove that these densities cross exactly once on $]0,m_\b[.$ Evaluating the densities and reparametrizing, this amounts to show that the function $x\,\mapsto\, c_1(1- c_2 x^\delta)^{1/\delta} + x - 1$ crosses the positive axis only once on $(0,1)$ for all $c_1, c_2 < 1 < \delta.$ Computing its second derivative $$-c_1c_2(\delta -1) x^{\delta-2}(1-c_2x^\delta)^{1/\delta -2}\, <\, 0$$ 
shows that the latter function is concave, negative at zero and positive at one, which finishes the proof.

\qed

\begin{Rem}\label{Rem3}

{\em (a) The decreasing property holds actually on $(0,\infty)$ for all $x \ge 0,$ and this is very easy to prove. Indeed, by definition one has
$$E_\alpha(\Ga(1+\a) x)\, =\, 1\, +\, x\, +\, \sum_{n\ge 2} \lpa \frac{(\Ga(1+\a))^n}{\Ga(1+n\a)}\rpa x^n$$
and for every $n\ge 2$ the mapping 
$$\a\,\mapsto \, \frac{(\Ga(1+\a))^n}{\Ga(1+n\a)}$$
has logarithmic derivative $n(\psi(1+\a) - \psi (1+n\a)) < 0$ for all $\a > 0$ and hence decreases on $(0,\infty).$                                                                                                                                                                                                                             The above result extends the property for $x\le 0$ and $\a\in (0,1),$ which is less trivial. Observe that the decreasing property does not hold anymore for $\a\in(1,\infty)$ because the Mittag-Leffler function has negative zeroes. For instance, the functions $E_1(-x) = e^{-x}$ and $E_2(-2x) = \cos (\sqrt{2x})$ are not comparable. Observe also that contrary to the functions $\cL_\a(x)$ - see again Proposition 10 in \cite{BS}, the non-rescaled functions $\Ea(x)$ are not comparable either, since $\Ea(0) = 1$ and $\a\mapsto \Ea'(0)$ is not monotonic on $(0,1).$  

\medskip

(b) It was shown in \cite{Frechet} that 
$$\frac{\M_\a}{\Ga(1-\a)}\,\prec_{st}\, \L\, \sim\, \lim_{\b\to 0} \lpa\frac{\M_\b}{\Ga(1-\b)}\rpa\quad\mbox{in distribution}$$
for all $\a\in (0,1),$ where $\L$ is the unit exponential random variable. We believe that the mapping
$$\a\,\mapsto\,\frac{\M_\a}{\Ga(1-\a)}$$
is non-increasing for the stochastic order on $(0,1)$, which would lead to 
$$\frac{1}{1+x}\, \le\,\Ea \lpa -\frac{x}{\Ga(1-\a)}\rpa \, \le\, E_\b \lpa -\frac{x}{\Ga(1-\b)}\rpa\, \le\, 1$$
for all $x\ge 0$ and $0<\a<\b <1,$ and hence refine the lower bound in \eqref{Unif}. It is easy to deduce from Theorem C that the mapping
$$\a\,\mapsto\,\Ea \lpa \frac{x}{\Ga(1-\a)}\rpa$$
decreases on $(0,1)$ for all $x > 0,$ and one can show by a log-convexity argument that the mapping
$$ \a\,\mapsto\, \esp\lcr \lpa \frac{\M_\a}{\Ga(1-\a)}\rpa^s\rcr\qquad\mbox{resp.} \qquad \a\,\mapsto\, \esp\lcr \lpa \frac{\M_\a}{\Ga(1-\a)}\rpa^{-s}\rcr$$
decreases on $(0,1)$ for all $s > 0$ resp. increases on $(0,1)$ for all $s\in (0,1).$ This is all in accordance with the above conjectured stochastic ordering property, which however still eludes us.

\medskip

(c) For all $d\in (0,2)$ and $p< 1/2,$ the random variable
$${\hat L}_1\, = \, \esp[{\hat L}_1]\,\prod_{k=0}^\infty \lpa\frac{k+1}{k+\a}\rpa \B_{\delta(1+\frac{k}{\a}),\delta(\frac{1}{\a}-1)}$$
with $\a =1-d/2\in (0,1)$ and $\delta =1-2p > 0$ is the (non-normalized) local time at time one and level zero of a noise-reinforced Bessel process with dimension $d$ and reinforcement parameter $p.$ See section 3 in \cite{Sun}, and \cite{Ber} for a whole study on the underlying local time process. If we consider the renormalized local time
$${\tilde L}_1(d,p)\, =\, \frac{{\hat L}_1}{\esp[{\hat L}_1]}\, \elaw\, \prod_{k=0}^\infty \lpa\frac{k+1}{k+\a}\rpa \B_{\delta(1+\frac{k}{\a}),\delta(\frac{1}{\a}-1)}\,\elaw\, \prod_{k=0}^\infty \lpa\frac{\Y_k}{\esp[\Y_k]}\rpa$$
with $\Y_k = \B_{\a(1+\frac{k}{\delta}),1-\a}^{\frac{\a}{\delta}}$ for all $k\ge 0$ - see (3.3) in \cite{Sun} for the second identity in law, the above proof shows that for each $p\in (-\infty,1/2),$ the mapping $d\mapsto {\tilde L}_1(d,p)$ is increasing for the convex order on $(0,2)$, with
$${\tilde L}_1(d,p)\;\claw\; 1\quad\mbox{as $d\to 0$}\quad\quad\mbox{and}\quad\quad {\tilde L}_1(d,p)\;\claw\; (1-2p)^{-1}\, \G_{1-2p}\quad\mbox{as $d\to 2$}$$
for each fixed $p\in (-\infty,1/2),$ where $\G_{1-2p}$ is the standard Gamma random variable. It is also easy to see from Theorem 1.2 in \cite{Ber} combined with the method of moments and Stirling's formula that for each fixed $d\in (0,2)$ one has
$${\tilde L}_1(d,p)\;\claw\; 1\quad\mbox{as $p\to-\infty$}\quad\quad\mbox{and}\quad\quad {\tilde L}_1(d,p)\;\claw\; c_\a^{-1}\, \cB(c_\a)\quad\mbox{as $p\to 1/2$}$$
where $c_\a = \a^{\a}/\Ga(1+\a) < 1$ by Gautschi's inequality on the Gamma function, and $\cB(q)$ stands for a Bernoulli random variable with parameter $q\in (0,1).$ An argument analogous to the above proof, and actually simpler, shows that for each $d\in (0,2),$ the mapping $p\mapsto {\tilde L}_1(d,p)$ is also increasing for the convex order on $(-\infty,1/2).$ We omit details.} 
\end{Rem}

\section{Proof of Theorem D} 

Throughout, we will use the usual Pochhammer notation $(a)_s = \Ga(a+s)/\Ga(a)$ for all $a > 0$ and $s > -a.$ In particular, we will sometimes set $(\b)_\a = \Ga(\a+\b)/\Ga(\b)$ for concision. Our argument relies on the fact that for all $\a\in (0,1], \b\ge\a$ and $x\in\rl$ one has
$$\Ga(\b)\, E_{\a,\b} \lpa \frac{\Ga(\a+\b)}{\Ga(\b)}\, x\rpa\, =\, \esp\lcr e^{x\,(\b)_\a \M_{\a,\b-\a}}\rcr$$
where $\M_{\a,\b-\a}$ is a positive random variable for which we refer to the introduction of \cite{FS2} for further details - see especially equations (2) and (3) therein. Setting 
$${\tilde\M_{\a,\b}} \, =\, (\b)_\a\, \M_{\a,\b-\a},$$
it follows from (4) in \cite{FS2} that $\esp[{\tilde\M}_{\a,\b}] = 1$ for all $\a\in (0,1]$ and $\b\ge\a.$ Reasoning as in Theorem C, we need to show that the mapping $\b\mapsto {\tilde\M}_{\a,\b}$ is non-decreasing for the convex ordering. To do so, we will use the independent multiplicative factorisation
$${\tilde\M}_{\a,\b} \, \elaw\, (\b)_\a\,\B_{\a,\b-\a}^\a\,\times\, \M_{\a,0}$$
for all $\a\in (0,1]$ and $\b >\a$, which is a direct consequence of (4) in \cite{FS2}, the known fact that
$$\esp\lcr\B_{\a,\b-\a}^{\a s}\rcr\, =\, \frac{(\a)_{\a s}}{(\b)_{\a s}}$$  
for all $s > -1,$ and Mellin inversion. By Corollary 3.A.22 in \cite{Shaked}, we are hence reduced to show that the mapping $\b\mapsto (\b)_\a\,\B_{\a,\b-\a}^\a$ increases for the convex ordering. Since
$$\esp\lcr (\b)_\a\,\B_{\a,\b-\a}^\a\rcr \, =\, \frac{\Ga(2\a)}{\Ga(\a)}$$
does not depend on $\b,$ we will use again the double intersection property given in (3.A.57) of \cite{Shaked}. The argument is the same as in the proof of Theorem C, but it is unfortunately longer because the involved densities exhibit more dependence on $\b.$  Fixing $\a\in (0,1]$ and setting $\rho^{}_{\b}$ for the density of the random variable $(\b)_\a\,\B_{\a,\b-\a}^\a,$ we first compute
$$\rho^{}_\b (x)\, =\, \kappa_{\b}\lpa 1 - \lpa\frac{x}{(\b)_\a}\rpa^{1/\a}\rpa^{\b-\a-1}\!\!\!\Un_{[0,\, (\b)_\a)}(x)$$
with
$$\kappa_{\b}\, =\, \frac{\Ga(\b)^2}{\Ga(\a+1)\Ga(\b-\a)\Ga(\a+\b)}\cdot$$
The mapping $\b\mapsto\kappa_{\b}$ increases since its log-derivative equals
$$2\psi(\b) \, -\, \psi(\b-\a)\, -\, \psi(\a+\b)\, =\, 2\a^2 \sum_{n\ge 0} \frac{1}{(n+\b) ((n+\b)^2 -\a^2)}\, >\, 0,$$
whereas the mapping $\b\mapsto(\b)_\a$ also increases by the log-convexity of the Gamma function. Therefore, if we fix $\b_2 > \b_1 > \a,$ we obtain
$$\frac{\rho^{}_{\b_1}(0)}{\rho^{}_{\b_2}(0)}\, =\, \frac{\kappa_{\b_1}}{\kappa_{\b_2}}\, <\, 1,\qquad\mbox{while}\qquad \rho^{}_{\b_1}(x) = 0\;\;\mbox{and}\;\; \rho^{}_{\b_2}(x) > 0\quad\mbox{for all $x\in [(\b_1)_\a,(\b_2)_\a),$}$$
and we are hence reduced to show that the strict sign sequence of the function $x\mapsto \rho^{}_{\b_2}(x) - \rho^{}_{\b_1}(x)$ is either $+-$ or $+-+$ on $(0,(\b_1)_\a).$ We will split the proof into four cases.

\begin{itemize}

\item Case $\b_2\ge\a+1 > \b_1.$ Here, $\rho^{}_{\b_2}$ is non-increasing on $(0,(\b_1)_\a),$ whereas $\rho^{}_{\b_1}$ increases towards $\infty,$ so that the sign sequence is $+-.$ 

\item Case $\b_2 > \a+1 = \b_1.$ Here, $\rho^{}_{\b_1}\equiv\kappa_{\a+1}$ on $(0,(\b_1)_\a),$ whereas $\rho^{}_{\b_2}$ decreases towards some positive constant which must be strictly smaller than $\kappa_{\a+1},$ since otherwise $\rho^{}_{\b_2}$ would not be a density. This shows that the sign sequence is again $+-.$  

\item Case $\b_2 > \b_1 > \a+1.$  Here, $\rho^{}_{\b_1}$ decreases to zero on $(0,(\b_1)_\a),$ whereas $\rho^{}_{\b_2}$ decreases towards some positive constant. We show that the equation $\rho^{}_{\b_1}(x) =\rho^{}_{\b_2}(x)$ has exactly two solutions on $(0,(\b_1)_\a),$ which implies that the sign sequence is $+-+.$ Observe that again, there must be at least two solutions since otherwise $\rho^{}_{\b_2}$ would not be a density. Setting $y = 1- (x/(\b)_\a)^{1/\a},$ the number of solutions amounts to that of
\begin{equation}
\label{SolCxcv}
1\, - \, \frac{(\b_1)_\a^{1/\a}}{(\b_2)_\a^{1/\a}}\, (1-y) \, -\,\kappa\, y^\mu\, =\, 0
\end{equation}
on $(0,1),$ with $\kappa = \lpa\kappa_{\b_1}/\kappa_{\b_2}\rpa^{\frac{1}{\b_2 - \a - 1}}$ and
$$\mu\, =\, \frac{\b_1-1-\a}{\b_2-1-\a}\, \in\, (0,1).$$
The function in \eqref{SolCxcv} is convex and there are hence at most two solutions to the equation.

\item Case $\a+1 > \b_2 > \b_1.$ Here, $\rho^{}_{\b_1}$ increases to $\infty$ on $(0,(\b_1)_\a),$ whereas $\rho^{}_{\b_2}$ increases towards some positive constant. Hence, the equation $\rho^{}_{\b_1}(x) =\rho^{}_{\b_2}(x)$ has at least one solution on $(0,(\b_1)_\a)$ and we show that there is only one, which implies that the sign sequence is $+-.$ We again need to solve \eqref{SolCxcv} but here one has $\mu > 1$ so that the function is concave, positive at zero and taking value $1-\kappa < 0$ at one, so it vanishes only once on $(0,1).$ 
 
\end{itemize}

This completes the proof of the increasing character of the mapping
$$\b\,\mapsto\, \Ga(\b)\, E_{\a,\b} \lpa \frac{\Ga(\a+\b)}{\Ga(\b)}\, x\rpa$$
on $(\a,\infty),$ and we finally compute its limit at $\infty$. If $x\ge 1,$ one has
$$\Ga(\b)\, E_{\a,\b} \lpa \frac{\Ga(\a+\b)}{\Ga(\b)}\, x\rpa \, \ge \, \Ga(\b)\, E_{\a,\b} \lpa \frac{\Ga(\a+\b)}{\Ga(\b)}\rpa\, = \, 2 \, +\, \sum_{n\ge 2} \frac{\Ga(\a+\b)^n}{\Ga(\b)^{n-1} \Ga(\b + n\a)}\; \to\; \infty$$
as $\b\to\infty,$ since $(\b)_\a\!\sim \b^\a$ and each summand tends to one. If $ x\le 0,$ Proposition 4 in \cite{BS} yields the bounds
\begin{equation}
\label{KSB}
\frac{1}{1 -\lpa\frac{\Ga(\a+\b)\Ga(\b-\a)}{\Ga(\b)^2}\rpa x}\; \le\; \Ga(\b)\, E_{\a,\b} \lpa \frac{\Ga(\a+\b)}{\Ga(\b)}\, x\rpa\; \le\; \frac{1}{1-x}
\end{equation}
and the convex ordering argument therein shows that these bounds extend to all $x < 1.$ This completes the proof since the lower bound of \eqref{KSB} converges to $1/(1-x)$ when $\b\to\infty.$

\qed

Theorem D implies the bounds
\begin{equation}
\label{Bd1}
\Ga(\a)\, E_{\a,\a} \lpa \frac{\Ga(2\a)}{\Ga(\a)}\, x\rpa \, \le\,\Ga(\b)\, E_{\a,\b} \lpa \frac{\Ga(\a+\b)}{\Ga(\b)}\, x\rpa\, \le\, \frac{1}{(1-x)_+}
\end{equation}
which are valid for all $\a\in (0,1], \b\ge \a$ and $x\in\rl,$ with the implicit notation that the upper bound is $\infty$ for $x\ge 1.$ The lower bound complements that of Proposition 4 in \cite{BS}, which converges to zero as $\b\to\a$ whereas here it converges to a non-trivial positive limit. It is natural to ask how the lower bound in \eqref{Bd1} behaves with respect to $\a$ and we have the following result, which is shown as Theorem C but will give some details.

\begin{Propo}
\label{Eaa}
For every $x\in\rl^\ast,$ the mapping 
$$\a\,\mapsto\, \Ga(\a)\, E_{\a,\a} \lpa \frac{\Ga(2\a)}{\Ga(\a)}\, x\rpa$$
decreases on $(0,1).$  
\end{Propo}

\begin{proof} With the notation used during the proof of Proposition \ref{Eab}, we have
$$\Ga(\a)\, E_{\a,\a} \lpa \frac{\Ga(2\a)}{\Ga(\a)}\, x\rpa\, =\, \esp\lcr e^{x (\a)_\a \M_{\a,0}}\rcr$$
for all $x\in\rl.$ Moreover, it follows from the last identity in law p.327 in \cite{FS2} that
$$(\a)_\a \M_{\a,0}\,\sim\, \prod_{n=0}^\infty \lpa \frac{n+1+\a}{n+2\a}\rpa \B_{2+\frac{n}{\a}, \frac{1}{\a}-1}\,\sim\,\prod_{n=0}^\infty \lpa \frac{\Ga(1+(n+2)\a)\Ga((n+2)\a)}{\Ga(1+(n+1)\a)\,\Ga((n+3)\a)}\rpa \B_{\a(n+2), 1-\a}^{\a}$$
where the second identity is obtained by a size-bias analysis as for (3.3) and (3.4) in \cite{Sun}. Each term in the product on the right-hand side is distributed as
$$\lpa\frac{n+3}{n+2}\rpa \T^{1,n+1}_\a$$
with the notation used during the proof of Theorem C, where it is shown that $\a\mapsto \T^{1,n+1}_\a$ decreases for the convex order. This shows that the mapping $\a\mapsto (\a)_\a \M_{\a,0}$ also decreases for the convex order, which concludes the proof as in Theorem C. 

\end{proof}

\begin{Rem}
\label{Rem4}
{\em (a) The above result implies the uniform bounds
$$e^{x}\, \le\, \Ga(\a)\, E_{\a,\a} \lpa \frac{\Ga(2\a)}{\Ga(\a)}\, x\rpa \, \le\,\Ga(\b)\, E_{\a,\b} \lpa  \frac{\Ga(\a+\b)}{\Ga(\b)}\, x\rpa\, \le\, \frac{1}{(1-x)_+}$$
for all $\a\in (0,1), \b \ge \a$ and $x\in\rl,$ which complete both those of \eqref{Bd1} and \eqref{Unif1}. Recall that
$$\Ga(\a)\, E_{\a,\a} \lpa \frac{\Ga(2\a)}{\Ga(\a)}\, x\rpa\,\to\, e^x\qquad \mbox{as $\a\to 1$}$$
and
$$\Ga(\b)\, E_{\a,\b} \lpa  \frac{\Ga(\a+\b)}{\Ga(\b)}\, x\rpa\, \to\, \frac{1}{(1-x)_+}
\qquad \mbox{as $\b\to \infty,$}$$
so that these bounds are optimal. On the other hand, it is easily shown that 
\begin{equation}
\label{LimitML}
\Ga(\a)\, E_{\a,\a} \lpa \frac{\Ga(2\a)}{\Ga(\a)}\, x\rpa\,\to\, \frac{1}{\lpa 1- \frac{x}{2}\rpa_+^2}\;\; \mbox{as $\a\to 0,$}
\end{equation}
which implies the better uniform bounds
\begin{equation}
\label{Bind}
e^x\, \le\, \Ga(\a)\, E_{\a,\a} \lpa \frac{\Ga(2\a)}{\Ga(\a)}\, x\rpa \, \le\, \frac{1}{\lpa 1- \frac{x}{2}\rpa_+^2}\cdot
\end{equation}
The above upper bound was already shown in Proposition 4 of \cite{BS} by other methods, together with the uniform lower bound
\begin{equation}
\label{Bind2}
\Ga(\a)\, E_{\a,\a} \lpa \frac{\Ga(2\a)}{\Ga(\a)}\, x\rpa \, \ge\, \frac{1}{\lpa 1- \sqrt{\frac{\Ga(1-\a)}{\Ga(1+\a)}} \frac{\Ga(2\a)}{\Ga(\a)}\, x\rpa^2}
\end{equation}
which holds for $x\le 0$ only. This bound is better than the above lower bound for large negative $x$, but worse for small negative $x$ or $\a$ close to 1. 

\medskip

(b) As mentioned in the introduction, for every $\a\in (0,1)$ and $\b\ge\a,$ the function $x\mapsto \Ga(\b) E_{\a,\b}(x)$ is a moment generating function and hence absolutely monotonic on $\rl.$ In particular, it has an inverse function which defines a bijection from $(0,\infty)$ onto $\rl,$ and can hence be viewed as a generalized logarithm. We refer to \cite{FS2} for conditions ensuring that this inverse function, denoted by $\log_{\a,\b},$ is a true logarithm, that is the function $x\mapsto -x\log_{\a,\b}(x)$ satisfies the three Shannon-Khinchin principles. A combination of Theorem C and Proposition \ref{Eaa} yields immediately the following uniform bounds
$$\frac{x}{x+1}\,\le\, \frac{\Ga(\b)}{\Ga(\a+\b)}\log_{\a,\b} (1+x)\,\le\,\frac{\Ga(\a)}{\Ga(2\a)}\log_{\a,\a} (1+x) \,\le\,\log(1+x)$$
for all $\a\in (0,1), \b \ge \a$ and $x > -1.$ Notice that the above limits for the corresponding Mittag-Leffler functions imply
$$\frac{\Ga(\a)}{\Ga(2\a)}\log_{\a,\a} (1+x) \,\to\,\log(1+x)\;\; \mbox{as $\a\to 1$}\quad\mbox{and} \quad\frac{\Ga(\b)}{\Ga(\a+\b)}\log_{\a,\b} (1+x)\, \to\, \frac{x}{x+1}
\;\; \mbox{as $\b\to \infty,$}$$
so that these uniform bounds are optimal. We also obtain from \eqref{LimitML}
$$\frac{\Ga(\a)}{\Ga(2\a)}\log_{\a,\a} (1+x) \,\to\,2\lpa 1 - \frac{1}{\sqrt{1+x}}\rpa\;\; \mbox{as $\a\to 0,$}$$
whence the uniform bounds 
$$\frac{x}{x+1}\, \le \, 2\lpa 1 - \frac{1}{\sqrt{1+x}}\rpa\,\le\,\frac{\Ga(\a)}{\Ga(2\a)}\log_{\a,\a} (1+x) \,\le\,\log(1+x)$$
for all $\a \in (0,1)$ and $x > -1.$}
\end{Rem}

We would like to conclude the paper with the following monotonicity result on the Mittag-Leffler $E_{\a,\a}(x)$ which is related to the main result of \cite{JZ}, and which also answers to some of the questions raised in Remark 3.2. therein.

\begin{Propo}
\label{Inv}
For every $\a\in (0,1),$ there exists $x_\a\in (0,1)$ such that the function $x\mapsto xE_{\a,\a} (-x)$ increases on $(0,x_\a)$ and decreases on $(x_\a,\infty).$ 
\end{Propo}

\begin{proof} Differentiating \eqref{ML Rep} and integrating by parts, we obtain
\begin{eqnarray*}
x^{\a}E_{\a,\a}(-x^\a) & = &  x\int_0^\infty e^{-xt}\lpa\frac{\sin(\pi\alpha)}{\pi(t^{\alpha}+ t^{-\a} +2\cos(\pi\alpha))} \rpa dt\\
& = & \a \int_0^\infty e^{-xt}\lpa\frac{\sin(\pi\alpha) (t^{-\a} - t^{\a})}{\pi \,t\,(t^{\alpha}+ t^{-\a} +2\cos(\pi\alpha))^2} \rpa dt, 
\end{eqnarray*}
which implies
\begin{equation}
\label{EaRep}
\frac{d}{dx} \lpa xE_{\a,\a} (-x) \rpa \, =\, x^{1/\a -1}\int_0^\infty e^{-x^{1/\a}t}\lpa\frac{\sin(\pi\alpha) (t^{\a} - t^{-\a})}{\pi (t^{\alpha}+ t^{-\a} +2\cos(\pi\alpha))^2} \rpa dt.
\end{equation}
On the other hand, the function 
$$t\,\mapsto\, \frac{\sin(\pi\alpha) (t^{\a} - t^{-\a})}{\pi (t^{\alpha}+ t^{-\a} +2\cos(\pi\alpha))^2}$$
has one strict sign-change (at $t=1$) on $(0,\infty),$ starting negative. This shows, again by Theorem 3.1 p.21 in \cite{Karlin}, that there exists $x_\a > 0$ with the required properties. Besides, we can deduce from the identity \eqref{EaRep}, as at the end of the proof of Theorem B, that
$$\frac{d}{dx} \lpa xE_{\a,\a} (-x) \rpa_{x=1} \, =\, \int_0^1 \lpa e^{-t} - t^{-2} e^{-1/t}\rpa\lpa\frac{\sin(\pi\alpha) (t^{\a} - t^{-\a})}{\pi (t^{\alpha}+ t^{-\a} +2\cos(\pi\alpha))^2} \rpa dt\; <\; 0,$$  
which implies $x_{\a} < 1.$

\end{proof}

Setting $m_\a = \sup_{x\ge 0} \lpa x E_{\a,\a} (-x)\rpa,$ we observe that \eqref{Bind} yields 
\begin{equation}
\label{BS} 
m_\a \; \le \; \frac{\Ga(2\a)}{2 \Ga(\a)^2} \; =\; \frac{4^{\a-1} \Ga(\a+1/2)}{\sqrt{\pi} \Ga(\a)}\; <\; \frac{\a}{2}
\end{equation}
where the second inequality follows from the log-convexity of the Gamma function. This upper bound improves on Proposition 3.1 in \cite{JZ}. Moreover, Proposition \ref{Inv} together with \eqref{Bind} and \eqref{Bind2} implies
the lower bound 
\begin{equation}
\label{BI}
m_\a \, \ge \, E_{\a,\a} (-1) \, \ge\, \lpa\frac{e^{-\frac{\Ga(\a)}{\Ga(2\a)}}}{\Ga(\a+1)}\vee\lpa\frac{1}{\sqrt{\Ga(1+\a)} +\sqrt{\Ga(1-\a)}}\rpa^2 \rpa\a.
\end{equation}
The function under brackets on the right-hand side has a positive minimum $\mu = 0.192744...$ and we hence have a lower bound $m_\a \ge \mu\a$ for all $\a\in (0,1).$ Observe finally that \eqref{BS} and \eqref{BI} give the sharp estimate $m_\a \sim\a/4$ as $\a\to 0.$

\section*{Acknowledgement}
This paper was partly written while the second author was visiting the University of Porto, whose support is gratefully acknowledged. The first author was supported by the ``Funda\c{c}\~{a}o para a Ci\^encia e a Tecnologia (FCT)" through the program ``Stimulus of Scientific Employment, Individual Support-6th Edition Call" with reference 2023.05967.CEECIND/CP2833/CT0001.

\end{document}